\documentclass[11pt]{amsart}
\usepackage{graphicx}
\usepackage{amssymb}
\usepackage{epstopdf}
\newtheorem{theorem}{Theorem}[section]
\newtheorem{proposition}[theorem]{Proposition}
\newtheorem{corollary}[theorem]{Corollary}

\newtheorem{question}[theorem]{Question}

\newtheorem{remark}[theorem]{Remark}
\def\primPsi{\Psi_{\text{\tiny prim}}}
\def\Mod{\mbox{Mod}}
\def\ModpA{\mbox{Mod}^{\text{\tiny pA}}}
\def\ModpAp{\mbox{Mod}^{\text{\tiny pA}+}}
\def\Msb{M_{\text{\tiny sb}}}

\def\sF{\mathcal{F}}

\def\R{\mathbb R}
\def\Z{\mathbb Z}

\def\sL{\mathcal L}
\def\sT{\mathcal T}
\def\sy{\mathcal{Y}}

\def\mod{\mathrm{mod}}
 \def\qed{\hfill\framebox(5,5){}}

\title[Small dilatation pseudo-Anosov mapping classes]{Small dilatation pseudo-Anosov mapping classes coming from the simplest hyperbolic braid.}
\author{Eriko Hironaka}
\begin{document}

\begin{abstract}  
 In this paper we study the minimum dilatation pseudo-Anosov mapping classes
arising from fibrations over the circle of a single 3-manifold, namely the mapping
torus for the "simplest hyperbolic braid". The dilatations that occur
include the minimum dilatations for orientable pseudo-Anosov mapping classes for genus
$g=2,3,4,5, \mbox{or} \ 8$.  In particular, we obtain the ``Lehmer example" in genus $g=5$,
and Lanneau and Thiffeault's conjectural minima in the orientable 
case for all genus $g$ satisfying $g= 2\ \mbox{or}\ 4 (\mod\ 6)$. Our examples show that the
minimum dilatation for orientable mapping classes is strictly greater than the
minimum dilatation for non-orientable ones when $g = 4,6, \mbox {or} \ 8$.   We also prove
that if $\delta_g$ is the minimum dilatation of pseudo-Anosov mapping classes on a genus
$g$ surface, then 
$$
\limsup_{g \rightarrow \infty} \ (\delta_g)^g \leq \frac{3 + \sqrt{5}}{2}.
$$
\end{abstract}

\maketitle

\section{Introduction}
Let $S_g$ be a closed oriented surface of genus $g \ge 1$, and let $\Mod_g$ be
the {\it mapping class group}, the group of orientation preserving homeomorphisms 
of $S_g$ to itself  up to isotopy.
A mapping class $\phi \in \Mod_g$ is called {\it pseudo-Anosov} if $S_g$ admits a pair of
$\phi$-invariant,  transverse measured, singular foliations
 on which $\phi$ acts by stretching transverse to one foliation
 by a constant $\lambda(\phi) > 1$ and contracting transverse to the other by $\lambda(\phi)^{-1}$.  The constant $\lambda(\phi)$ is called the {\it (geometric) dilatation} of $\phi$.   
 A mapping class is pseudo-Anosov if it is
neither periodic nor reducible  \cite{Thurston88, FLP79, CB88}.   
Denote by $\ModpA_g$ the set of pseudo-Anosov mapping classes in $\Mod_g$.

A pseudo-Anosov mapping class $\phi$ is defined to be {\it orientable} if its invariant
foliations are orientable.  We will denote the set of orientable pseudo-Anosov mapping
classes by $\ModpAp_g$.
Let $\lambda_{\mathrm{hom}}(\phi)$ be the
spectral radius of the action of $\phi$ on the first homology of $S$.
Then 
$$
 \lambda_{\mathrm{hom}}(\phi) \leq \lambda(\phi),
$$ 
with equality if and only if $\phi$ is orientable
(see, for example, \cite{LT09}(p. 5) and \cite{KS:Reps} (Theorem 1.4)).

The dilatations $\lambda(\phi)$, for $\phi \in \ModpA_g$ satisfy reciprocal monic integer polynomials of
degree bounded from above
by $6g-6$ \cite{Thurston88}.  If $\phi$ is orientable the degree is bounded by $2g$.
For fixed $g$, it follows that $\lambda(\phi)$ achieves a minimum $\delta_g > 1$ 
on $\ModpA_g$ (see also, \cite{AY,Ivanov90}).    Let $\delta_g^+$ be the minimum dilatation
among orientable  elements of $\ModpAp_g$.
 
In this paper, we address the following question (cf. \cite{Penner91,McMullen:Poly, Farb:Problems}):

\begin{question} What is the behavior of $\delta_g$ and $\delta_g^+$ as functions of $g$?
\end{question}

So far, exact values of $\delta_g$  have only been found for $g \leq 2$.
For $g=1$,  $\Mod_1 = \mathrm{SL}(2;\Z)$, and
$$
\delta_1 = \frac{3+\sqrt{5}}{2}.
$$
For $g=2$, Cho and Ham \cite{CH08} show that $\delta_2$ is the largest
real root of
$$
t^4 - t^3 - t^2 - t + 1 = 0
$$
($\approx 1.72208$).  

 In the orientable case more is known due to recent
results of Lanneau and Thiffeault \cite{LT09}.
Given  $(a,b) \in \Z \oplus \Z$ with $0 <  a < b$, let
$$
LT_{(a,b)}(t) = t^{2b} - t^b(1 + x^a + x^{-a}) + 1,
$$
and let $\lambda_{(a,b)}$ be the largest real root of 
$LT_{(a,b)}(t)$.   

\begin{theorem}[Lanneau-Thiffeault \cite{LT09} Theorems 1.2 and 1.3] \label{LT-thm}
For $g = 2,3,4,6$, and $8$,
$$
\lambda_{(1,g)} \leq \delta_g^+
$$
with equality when $g = 2,3$, or $4$.  
\end{theorem}

\noindent
For $g=2$, the value of $\delta_2^+$ was first determined by Zhirov \cite{Zhirov95}.
For $g=5$, Lanneau and Thiffeault show that $\delta_5^+$ equals Lehmer's
number ($\approx 1.17628$) \cite{Lehmer33}.  This dilatation is realized as a product of multi-twists
along a curve arrangement dual to the $E_{10}$ Coxeter graph (see \cite{Leininger04, Hironaka:Coxeter}),
and as the monodromy of the (-2,3,7)-pretzel knot (see \cite{Hironaka98}).   Lanneau and 
Thiffeault also find a lower bound for $\delta_7^+$.   An example realizing this bound
can be found in \cite{AD10} (p.4) and \cite{KT10} (Theorem 1.12).

Based on their results, Lanneau and Thiffeault ask:

\begin{question}[\cite{LT09} Question 6.1]\label{LT-ques}
Is $\delta_g^+ = \lambda_{(1,g)}$ for all even $g$?
\end{question}

For convenience, we will call the affirmative answer to their question the {\it LT-conjecture}.

In our first result, we improve on the following
previous  best bounds for the minimum dilatation of infinite families 
$$
(\delta_g)^g \leq (\delta_g^+)^g \leq 2 + \sqrt{3}
$$
found in \cite{Minakawa:dil, HK:braidbounds}.

\begin{theorem}\label{bounds-thm} If $g = 0,1,3,$ or $4 (\mod \ 6)$, $g \ge 3$,  then 
$$
\delta_g  \leq \lambda_{(3,g+1)},
$$
and if $g = 2,$ or $5(\mod \ 6)$ and $g \ge 5$, then
$$
\delta_g \leq \lambda_{(1,g+1)}.
$$
\end{theorem}

\noindent
For the orientable case, our results complement those of Lanneau and
Thiffeault for $g = 2,$ or $4(\mod \ 6)$.

\begin{theorem}\label{orientable-bounds-thm}
Let  $g \ge 3$.  Then
$$
\begin{array}{ll}
\delta_g^+ \leq \lambda_{(3,g+1)}&\quad\mbox{if}\ g=1,\ \mbox{or}\ 3 (\mod \ 6),\\
\delta_g^+ \leq \lambda_{(1,g)} &\quad\mbox{if}\ g=2,\ \mbox{or} \ 4(\mod \ 6),
\ \mbox{and}\\
\delta_g^+ \leq \lambda_{(1,g+1)} &\quad\mbox{if}\ g = 5 (\mod \ 6).
\end{array}
$$
\end{theorem}

\noindent

Putting Theorem~\ref{orientable-bounds-thm} together with Lanneau and Thiffeault's lower bound for
$g=8$ gives:

\begin{corollary}\label{Eight-cor} For $g=8$, we have
$$
\delta_8^+ = \lambda_{(1,8)}.
$$
\end{corollary}

The following is a table of the minimal dilatations that arise in this paper's  examples for
genus $1$ through $12$.   All numbers in the table are truncated  to 5 decimal places.  
An asterisk $*$ marks the numbers that
have been verified to equal $\delta_g^+$ (resp., $\delta_g$). 
For singularity-type, we use the convention that $(a_1,\dots,a_k)$ means
that the singularities of the invariant foliations have degrees $a_1,\dots,a_k$
(see Lanneau and Thiffeault's notation  \cite{LT09} p.3).   
The singularity-types 
for our examples are derived from the formula given in Proposition~\ref{prongs-prop}.
\bigskip

\begin{center}
\begin{tabular}{|c|| l| l| |l| l|}
\hline
$g $& orientable& degrees of singularities &unconstrained &degrees of singularities\\
\hline
1 & $2.61803{}^*$ & no sing. &$2.61803{}^*$& no sing.\\
\hline
2 & $1.72208{}^*$ & $(4)$ & $1.72208{}^*$ & $(4)$\\
\hline
3 & $1.40127{}^* $&$(2,2,2,2)$& 1.40127 &(2,2,2,2)\\
\hline
4 & $1.28064{}^*$ &(10,2)& 1.26123&(3,3,3,3) \\
\hline
5 & $1.17628 {}^*$&(16)& 1.17628&(16)\\
\hline
6 & - &-&1.1617&(5,5,5,5)\\
\hline
7 & 1.13694& (6,6,6,6)& 1.13694& (6,6,6,6)\\
\hline
8 & $1.12876{}^*$&(22,6)& 1.1135&(25,1,1,1)\\
\hline
9& 1.1054 &(8.8.8.8)& 1.1054 &(8,8,8,8)\\
\hline
10 & 1.10149 &(28,8)& 1.09466 &(9,9,9,9)\\
\hline
11 &1.08377 &(34,2,2,2)&1.08377&(34,2,2,2)\\
\hline
12 & -&-&1.07874&(11,11,11,11)\\
\hline
\end{tabular}

\bigskip
Table 1: Minimal orientable and unconstrained dilatations coming from $\Msb$ 
\end{center}
\bigskip

For $g=2,3,4,$ and $5$, our orientable examples agree both in dilatation and in singularity-type
with the previously known minimizing examples (see \cite{LT09} \S 3, \S 4, \S 6).     
For $g=8$, our example agrees with the singularity-type anticipated by Lanneau and Thiffeault
\cite{LT09} (6.4).  
We prove that the known minimial dilatation examples for $g=2,3,4,5,$ and $8$ arise as the monodromy
of fibrations of a single 3-manifold $\Msb$.
For $g=7$, our minimal example gives a strictly larger dilatation than $\delta_7^+$.
(The dilatation $\delta_7^+$ is realized in \cite{KT10} and  \cite{AD10}.)

Lanneau and Thiffeault show that $\delta_5^+ \le \delta_6^+$, and hence 
$\delta_g^+$ is not strictly monotone decreasing (cf. \cite{Farb:Problems} Question 7.2).
Theorem~\ref{orientable-bounds-thm} implies the following stronger statement.

\begin{proposition}
If the LT-conjecture is true, then $\delta_g^+ \leq \delta_{g+1}^+$, whenever $g = 5 (\mod \ 6)$.
\end{proposition}

Another consequence concerns the question of whether the inequality 
$\delta_g \leq \delta_g^+$ is strict for any or all $g$.   In \cite{KT10} and \cite{AD10} it is
shown that $\delta_5 < \delta_5^+$.  Table 1 shows
the following.

\begin{corollary}\label{inequality-cor} For $g=4,6$, and $8$ we have
$$
\delta_g < \delta_g^+.
$$
\end{corollary}

\noindent
If the LT-conjecture is true, then Theorem~\ref{bounds-thm} and
Proposition \ref{monotonicity-prop} imply that
the phenomena revealed in
Corollary~\ref{inequality-cor}
repeats itself periodically.

\begin{proposition}
If  the LT-conjecture is true, then for all even $g \ge 4$  we have
$$
\delta_g < \delta_g^+.
$$
\end{proposition}

For large $g$, it is known that $\delta_g$ and $\delta_g^+$
converges to $1$. Furthermore, 
\begin{eqnarray}\label{Penner-eqn}
\log(\delta_g) \asymp \frac{1}{g} \qquad \mbox{and}\qquad \log(\delta_g^+) \asymp \frac{1}{g}
\end{eqnarray}
(see \cite{Penner91, McMullen:Poly, Minakawa:dil, HK:braidbounds}).
The LT-conjecture together with (\ref{Penner-eqn}) leads to the natural question:

\begin{question} [\cite{McMullen:Poly} p.551, \cite{Farb:Problems} Problem 7.1]  Do the sequences
$$
(\delta_g)^g \qquad\mbox{and}\qquad (\delta_g^+)^g
$$
converge as $g$ grows?  What is the limit?
\end{question}

Theorem 1.3 and Theorem 1.4 imply the following.

\begin{theorem}\label{limit-thm}
$$
\limsup_{g \rightarrow \infty}\ (\delta_g)^g \leq  \frac{3+\sqrt{5}}{2} 
$$
and
$$
\limsup_{g \neq 0 (\mod 6)}\ (\delta_g^+)^g \leq \frac{3 + \sqrt{5}}{2}.
$$
\end{theorem}

This leads to the following question.

\begin{question}[Golden Mean Question]\label{GM-ques}  Do the sequences
$(\delta_g)^g$ and  $(\delta_g^+)^g$ satisfy
$$
\lim_{g \rightarrow \infty} (\delta_g)^g =
\lim_{g \rightarrow \infty} (\delta_g^+)^g =
\frac{3 + \sqrt{5}}{2} = 1 + \mbox{golden mean}\ ?
$$
\end{question}

For any pseudo-Anosov mapping class $\phi$, let $M(\phi)$ be the mapping torus of $\phi$. 
Conversely, given a compact hyperbolic 3-manifold with torus boundary components $M$,
let $\Phi(M)$ be the collection of pseudo-Anosov  mapping classes $\phi$
such that $M = M(\phi)$.

We prove Theorem~\ref{bounds-thm} and Theorem~\ref{orientable-bounds-thm}
by exhibiting a family of mapping classes
$\phi_{(a,b)} \in \Phi(\Msb)$  for a single hyperbolic 3-manifold $\Msb$.

 Let $\Sigma$ be
the suspensions of singularities of the stable and unstable foliations
of $\phi$ and let 
$$
M^*(\phi) = M(\phi) \setminus \Sigma.
$$

\begin{theorem}[\cite{FLM09} Theorem 1.1]\label{FLM-thm}
The set
$$
\sT_P = \{M^*(\phi) \ : \  \phi \in \ModpA_g, \lambda(\phi) \leq P^\frac{1}{g}\}
$$
is finite for any $P > 1$.
\end{theorem}

The
asymptotic equations (\ref{Penner-eqn}) and Theorem~\ref{FLM-thm} imply that
$$  
\sT = \{M^*(\phi) \ :\ \phi \in \ModpA_g,  \lambda(\phi) = \delta_g\} 
$$
and
$$
\sT^+ = \{M^*(\phi) \ : \ \phi \in \ModpAp_g \lambda(\phi) = \delta^+_g\} 
$$
are finite.

This invites the question:

\begin{question} How large are the sets $\sT$ and $\sT^+$?
\end{question}

If the LT-conjecture is true, then our results imply that a single 3-manifold would realize
$\delta_g^+$ for all $g = 2,4 (\mod\ 6)$.   The manifold we study in this paper
 $\Msb$ is the 
complement of the $6{}_2^2$ braid (see Rolfsen's tables \cite{Rolfsen76}, and
Figure~\ref{link-fig}).   Another 3-manifold that produces small dilatation
mapping classes is the complement
$M_{\mbox{-2,3,8}}$ of the $(-2,3,8)$-pretzel link in $S^3$.  These  have been
 studied independently by
Kin and Takasawa \cite{KT10}
and Aaber and Dunfield \cite{AD10}.   For certain genera the mapping classes in
$\Phi(M_{\mbox{-2,3,8}})$ have smaller dilatation than the minima
realized by $\Msb$, but the asymptotic behavior of the 
minimal dilatations for large genus, 
supports the affirmative to Question~\ref{GM-ques}.
Both $M_{-2,3,8}$ and $\Msb$ can be obtained from the 
{\it magic manifold} by Dehn fillings \cite{MP06}.   The pseudo-Anosov braid monodromies
with smallest known dilatations found in \cite{HK:braidbounds}
are also realized on the magic manifold 
\cite{KT:magicmanifold}.

Section~\ref{background-section} contains a brief review of
Thurston norms, fibered faces
and the Teichm\"uller polynomial.  These are the basic tools used in this paper.  In
Section~\ref{example-section} we describe our family of examples, and in
Section~\ref{growth-section} we
prove Theorem~\ref{bounds-thm} and Theorem~\ref{orientable-bounds-thm}.

\bigskip
\noindent
{\bf Acknowledgments:}
The author thanks Curt McMullen and Thomas Koberda for many helpful conversations,
and thanks Eiko Kin, Spencer Dowdall 
and the referee of this article for corrections to earlier versions.
The author is grateful for the hospitality of the Harvard Mathematics Department, where
she wrote this paper as a visiting scholar in Fall 2009.

\section{Background and tools}\label{background-section}

In this section we give a brief review of invariants and properties of fibrations of a hyperbolic3-manifold $M$, emphasizing the tools that we will use in the rest of the paper.
For more details see, for example,  \cite{Thurston86, FLP79, McMullen:Poly, McMullen:Alex}.

The theory of fibered faces of the Thurston norm ball and  the existence of
Teichm\"uller polynomials provides a way to study in a single picture a
collection of  pseudo-Anosov
mapping classes defined on surfaces of different Euler characteristics and
genera.  Assume $M$ is a compact hyperbolic 3-manifold with boundary.
  Given an embedded surface $S$ on $M$, let
 $\chi_- (S)$ be the sum of $|\chi(S_i)|$, where $S_i$ are the connected components of $S$
 with negative Euler characteristic.
 The  Thurston norm of $\psi \in \mathrm{H}^1(M;\Z)$ is defined to be 
 $$
 || \psi ||_T = \min \chi_-(S),
 $$
 where the minimum is taken over oriented embedded surfaces
  $(S, \partial S) \subset (M,\partial M)$ such that the
  class of $S$ in $ \mathrm{H}_2(M, \partial M; \Z)$ is dual to $\psi$.
  
Elements of $\mathrm{H}^1(M;\Z)$ are canonically associated with
epimorphisms
 $$
 \pi_1(M;\Z) \rightarrow \Z.
$$
We thus make the following natural identification:
 $$
 \mathrm{H}^1(M;\Z) = \mathrm{Hom}(\pi_1(M), \Z) = \mathrm{Hom}(\mathrm{H}_1(M;\Z) ,\Z).
 $$
 This can be considered as a lattice $\Lambda_M$ inside $\R^{b_1(M)}$.
 If $\psi \in \Lambda_M$ corresponds to a fibration
 $$
\psi: M \rightarrow S^1
 $$
 we say that $\psi$ is {\it fibered}.   In this case the Thurston norm of $\psi$ is given by
 $$
 ||\psi||_T = \chi_-(S),
 $$
 where $S$ is homeomorphic to the fiber of $\psi$.
 Let $$
 \Psi(M) = \{ \psi : M \rightarrow S^1\ : \ \mbox{$\psi$ is a fibration}\}.
 $$
The {\it monodromy} $\phi$ of $\psi \in \Psi(M)$ is the mapping class $\phi : S \rightarrow S$,
such that  $M$ is the mapping torus of $\phi$, and $\psi$ is the natural projection
to $S^1$.   Since $M$ is hyperbolic,  $\phi$ is automatically pseudo-Anosov.
  
Let $B$ be the
 unit ball in $\R^{b_1(M)}$ with respect to the extended Thurston norm.
 
 \begin{theorem}[\cite{Thurston86}]  The Thurston norm ball
 $B$ is a convex polyhedron and for any top-dimensional open face
 $F$ of $B$,
$(F \cdot \R^+) \cap \Psi(M)$ is either empty or equal to $(F\cdot \R^+) \cap \Lambda_M$.
\end{theorem}

If $(F \cdot \R^+) \cap \Psi(M) \neq \emptyset$, we say $F$ is a {\it fibered face} of $B$.
An element of $\Psi(M)$ is called {\it primitive} if its fiber is connected.
The elements of $\Lambda_M$ project to the  rational points on the boundary of
$B$.  If $F$ is a fibered
 face, then each rational point $x$ on $F$ corresponds to a unique primitive 
 element $\psi_x \in \Psi(M)$, namely the element of $(x\cdot \R^+) \cap \Psi(M)$ that
 lies closest to the origin.

\begin{theorem}[\cite{Fried82}, Theorem E]\label{Fried-thm} There is a continuous function
$\sy$, homogeneous of degree one, defined on the  fibered cone in $\R^{b_1(M)}$, so that if $\psi$ is fibered with monodromy $\phi_\psi$, then
$$
\sy(\psi) = \frac{1}{\log(\lambda(\phi_\psi))}.
$$
The function $\sy$ is concave and tends to zero along the boundary of the cone.
\end{theorem}

\begin{corollary} For each fibered face $F$,
$$
\overline \lambda (\psi) = \lambda(\phi_\psi)^{||\psi||_T},
$$
extends to a continuous function on $F \cdot \R^+$  that is constant
on rays through the origin, and $\overline\lambda$
achieves a unique minimum on $F$.
\end{corollary}

Let $G$ be a group and $\psi : G \rightarrow \Z$ a homomorphism.  If $f \in \Z[G]$,
$$
f = \sum_{g \in G} \alpha_g g,
$$
then the {\it specialization} of $f$ at $\psi$ is the polynomial in $\Z[t]$ defined by
$$
f^\psi(t) = \sum_{g \in G} \alpha_g \psi(g),
$$
where we think of $\psi(g) \in \Z = \langle t \rangle$ as a power of $t$.

For a monic integer polynomial $p(x)$, the {\it house} of $p(x)$, written $h(p)$, is the
absolute value of the largest root of $p$.

\begin{theorem}[\cite{McMullen:Poly}]  Let $F$ be a fibered face for a 3-manifold
$M$, and let $G = H_1(M;\Z)$.  Then there is an element $\theta_F \in \Z[G]$ such
that for all integral lattice points $\psi$ in the fibered cone of $F$,
$$
\lambda(\phi_\psi) = h(\theta_F^\psi).
$$
\end{theorem}

The polynomial $\theta_F$ is called the {\it Teichm\"uller polynomial} of $M$
for the fibered face $F$.

 \section{The mapping torus for the simplest hyperbolic braid}\label{example-section}
 
We now look at a particular 3-manifold, and study properties of its fibrations.  This example
has also been studied in (\cite{McMullen:Poly} \S 11), and the first part of this section
will be a review of what is found there.
\begin{figure}[htbp]
\begin{center}
\includegraphics[height=5cm]{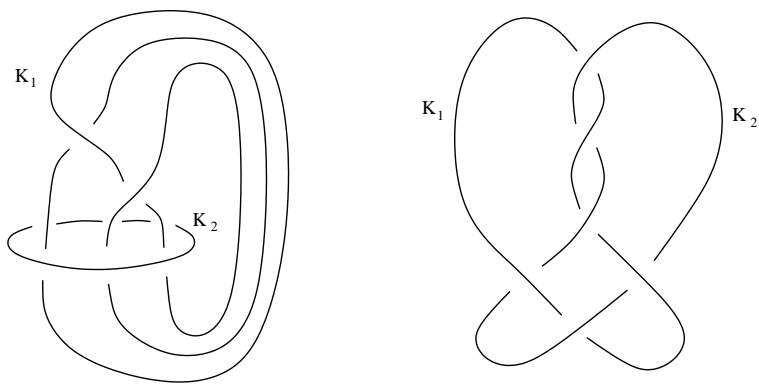}
\caption{Two diagrams for the link $6{}^2_2$.}
\label{link-fig}
\end{center}
\end{figure}

\begin{figure}[htbp]
\begin{center}
\includegraphics[height=5cm]{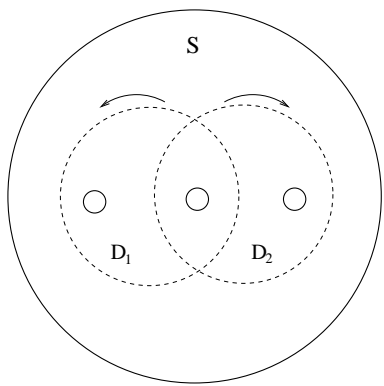}
\caption{Braid monodromy associated to $\sigma_1 \sigma_2^{-1}$.} 
\label{braidmonodromy-fig}
\end{center}
\end{figure}

Let $M = S^3 \setminus N(L)$, where
$L$ is the link drawn in two ways in Figure~\ref{link-fig}, and $N(L)$ is a tubular neighborhood.
As seen from the left diagram in Figure~\ref{link-fig},  $M$
 fibers over the circle with fiber a sphere with four boundary components $S_{0,4}$.  
 Let $\psi_0 : M \rightarrow S^1$ be the corresponding fibration, and
let $\phi_0 : S_{0,4} \rightarrow S_{0,4}$ be the monodromy.  Then $\phi_0$
is the mapping class associated to the braid written with respect to standard generators as
$\sigma_1 \sigma_2^{-1}$ (see Figure~\ref{braidmonodromy-fig}) and its dilatation
is given by
$$
\lambda(\phi_0) = \frac{3 + \sqrt{5}}{2}.
$$
The braid $\sigma_1 \sigma_2^{-1}$  has also been called the ``simplest pseudo-Anosov braid"
(see \cite{McMullen:Poly} \S 11).

Let $K_1$ and $K_2$ 
be the components of $L$ as drawn in Figure~\ref{link-fig}.   Let
$\mu_1$ be the meridian of $K_1$ and $\mu_2$ be the meridian of $K_2$.
Then any element $\psi$ of $H^1(M;\Z)$ is determined by its values 
$$
(a,b) = (\psi(\mu_1),\psi(\mu_2)) \in \Z \times \Z.
$$

With respect to these coordinates, the Thurston norm and the Alexander norm both are given by
\begin{eqnarray}\label{Thurston-norm-eqn}
||(a,b)|| = \max\{ 2|a|, 2|b|\}.
\end{eqnarray}

\begin{figure}[htbp]
\begin{center}
\includegraphics[height=8cm]{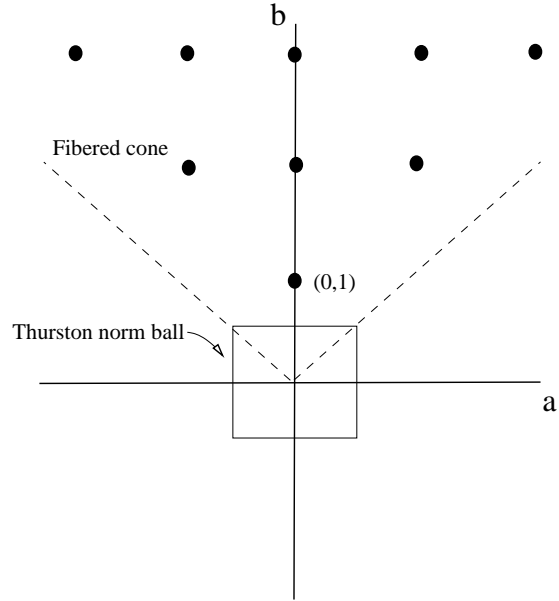}
\caption{Fibered cone $\Psi$ containing $\psi = (0,1)$.}
\label{FiberedFace-fig}
\end{center}
\end{figure}

The lattice points $\Lambda_M$ in the fibered cone $F \cdot \R^+$ 
defined by $\psi = (0,1)$ is the set
$$
\Psi  = \{(a,b) \in \Z \times \Z \ : \ b > 0,  -b < a < b\}
$$
as shown in Figure~\ref{FiberedFace-fig}.    For the rest of this paper, we
will only be concerned with the subset $\primPsi \subset \Psi$ consisting of
 elements of $\Psi$ with connected fibers, i.e., the {\it primitive elements}.
Thus,
$$
\primPsi = \{(a,b) \in \Z \times \Z \ : \ b > 0, -b < a < b, \gcd(a,b) = 1\}.
$$

The Alexander polynomial for $L$  is given by
\begin{eqnarray}\label{AlexPoly-eqn}
\Delta_L(x,u) = u^2 -  u( 1 -  x -  x^{-1})+ 1
\end{eqnarray}
(see Rolfsen's table  \cite{Rolfsen76}), and the Teichm\"uller polynomial is given by
\begin{eqnarray}\label{Teich-eqn}
\Theta_L(x,u) = u^2 - u(1 + x + x^{-1})  + 1
\end{eqnarray}
(see \cite{McMullen:Poly} p.47).

Specialization to the element $(a,b) \in \mathrm{H}^1(M;\Z)$ is the same as plugging $(t^a,t^b)$ into the equations for the Alexander and
Teichm\"uller polynomials (see Section~\ref{background-section}).

\begin{proposition}\label{dilatation-prop} If $(a,b) \in \primPsi$, then the associated monodromy
$\phi_{(a,b)}$ is pseudo-Anosov with geometric dilatation equal to the largest
root $\lambda_{(a,b)}$ of 
$$
\Theta_L(t^a,t^b) = t^{2b} - t^b(1 + t^a + t^{-a}) + 1,
$$
and homological dilatation the maximum norm
among roots of the polynomial
$$
\Delta_L(t^a,t^b) = t^{2b} - t^b(1 - t^a -t^{-a}) + 1.
$$
\end{proposition}

\begin{corollary}\label{orientability-cor} If $(a,b) \in  \primPsi$, then the associated
monodromy $\phi_{(a,b)}$ is orientable if $a$ is odd and $b$ is even.
\end{corollary}

\noindent
{\bf Proof.}  If $a$ is odd and $b$ is even, then 
the roots of $\Theta_L(t^a,t^b)$ are the negatives of the roots  of $\Delta_L(t^a,t^b)$.
This implies that the geometric and homological dilatations of $\phi_{(a,b)}$ are
equal, and therefore $\phi_{(a,b)}$ is orientable.
\qed

Later in this section, we prove the converse of Corollary~\ref{orientability-cor}.  First we consider
how the monodromy behaves near the boundary of $S_{(a,b)}$.

\begin{proposition}\label{boundary-prop}  
Let $\phi_{(a,b)} : S_{(a,b)} \rightarrow S_{(a,b)}$ be the monodromy associated to $(a,b) \in \primPsi$.
The boundary components
of $S_{(a,b)}$ has $\gcd(3,a)$ components coming from
$T(K_1)$ and $\gcd(3,b)$  coming from $T(K_2)$.  Thus, the total number of boundary
components of $S_{(a,b)}$ is given by
$$
\left \{
\begin{array}{ll}
2 & \mbox{if $\gcd(3,ab) = 1$}\\
4 &\mbox{if $\gcd(3,ab)=3$}
\end{array}
\right .
$$
\end{proposition}

\noindent
{\bf Proof.}
The number of components in $T(K_i) \cap S_{(a,b)}$ is the
index of the image of $\pi_1(T(K_i))$ in $\Z$ under the composition of maps 
$$
\pi_1(T(K_i)) \rightarrow \pi_1(M) \rightarrow \Z
$$
induced by inclusion and $\psi_{(a,b)}$.

For $i=1,2$, let  $\ell_i$ be the longitude of $K_i$ that is contractible in $S^3 \setminus K_i$.
Then, for $T(K_1)$ we have 
$$
\psi_{(a,b)}(\mu_1) = a \quad \mbox{and}\quad \psi_{(a,b)}(\ell_1) = 3 \psi_{(a,b)}(\mu_2) = 3b,
$$
so the number of boundary components contributed by $T(K_1)$ is 
$$
\gcd (a,3b) = \gcd(3,a),
$$
since we are assuming that $\gcd(a,b) = 1$.
The
contribution of $T(K_2)$ is computed similarly.
\qed

\begin{proposition}\label{genus-prop}
The genus of $S_{(a,b)}$, for $(a,b) \in \primPsi$ is given by
\begin{eqnarray*}
g(S_{(a,b)}) &=& |b| +\left (1 - \frac{\gcd(3,a) + \gcd(3,b)}{2} \right )\\
&=&
\left \{
\begin{array}{ll}
|b|  &\qquad\mbox{if $\gcd(3,ab)=1$}\\
|b|-1 &\qquad\mbox{if $\gcd(3,ab)=3$}.\\
\end{array}
\right. 
\end{eqnarray*}
\end{proposition}

\noindent
{\bf Proof.} From (\ref{Thurston-norm-eqn}) we have
$$
2|b| = \chi_-(S_{(a,b)}) = 2g - 2 + \gcd(3,a) + \gcd(3,b).
$$
\qed

\begin{figure}[htbp]
\begin{center}
\includegraphics[height=5cm]{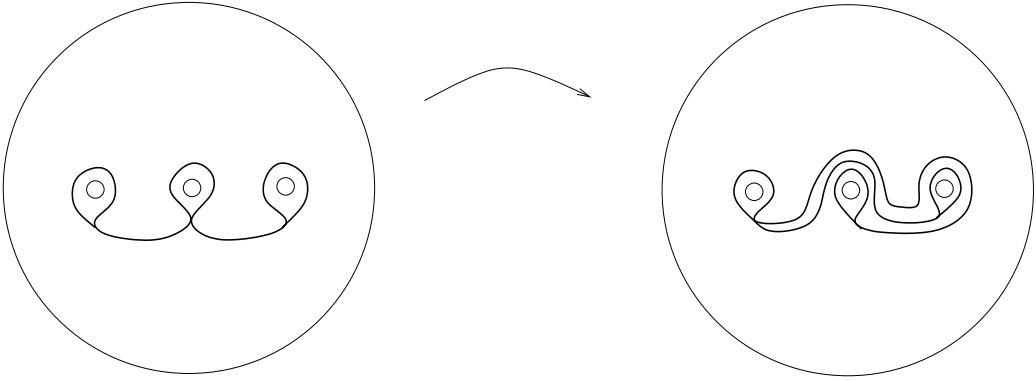}
\caption{Train track for $\phi: S \rightarrow S$.}
\label{traintrack-fig}
\end{center}
\end{figure}

\begin{proposition}\label{prongs-prop} Let  $(a,b) \in \primPsi$, and let $\sF$ be a
 $\phi_{(a,b)}$-invariant
foliation.  Then $\sF$ 
\begin{enumerate}
\item has no interior singularities,
\item is  $(3b/\gcd(3,a))$-pronged at each of the $\gcd(3,a)$ boundary components coming from $T(K_1)$, and
\item  is $(b/\gcd(3,b))$-pronged at each of the $\gcd(3,b)$ boundary components coming  from $T(K_2)$.
\end{enumerate}
\end{proposition}

\noindent
{\bf Proof.}  Let $\sL$ be the lamination of $M$ defined by suspending $\sF$ over
$M$ considered as the mapping torus of $\phi$.  From the train track for
$\phi$ (Figure~\ref{traintrack-fig}), one sees 
that each of the boundary components of $S$ are one-pronged, and that
there are no other singularities.   It follows that $\sL$ has no singularities
outside a neighborhood of the $K_i$, and near
each $K_i$ the leaves of $\sL$ come together at a simple closed curve $\gamma_i \in \mathrm{H}_1(T(K_i))$.  Write
$$
\gamma_i = r_i \mu_i +s_i \ell_i
$$
for $i=1,2$.

For $(a,b) \in \primPsi$, the number of intersections of $\gamma_i$ with $S_{(a,b)}$ is 
the image of $\gamma_i$ under the epimorphism
$$
\psi_{(a,b)} : \pi_1(M) \rightarrow \Z
$$
defining the fibration.
Figure~\ref{traintrack-fig} shows that
 $s_1 = 1$ and $r_2 = 1$.   Using the identities
\begin{eqnarray*}
s_1 = 1&\qquad&
\ell_1 = 3 \mu_2 ,\\
r_2 = 1&\qquad&
\ell_2 = 3 \mu_1,
\end{eqnarray*}
we have
\begin{eqnarray*}
\psi_{(a,b)}(\gamma_1) &=& r_1 \psi_n(\mu_1)  + 3 \psi_n(\mu_2) = r_1a +3b\\
\psi_{(a,b)}(\gamma_2) &=&  \psi_n(\mu_2) + 3 s_2 \psi_n(\mu_1) = 3s_2a + b.
\end{eqnarray*}
Let $m_1=\gcd(3,a)$ and $m_2 = \gcd(3,b)$.
Then $\phi_{(a,b)}$ is $(r_1a + 3b)/m_1$-pronged at
$m_1$ boundary components and  $(3s_2a +b)/m_2$-pronged
at $m_2$ boundary components.
We find $r_1$ and $s_2$ by looking at some particular examples.

In general, if $f: \Sigma \rightarrow \Sigma$ is pseudo-Anosov
on a compact oriented surface $\Sigma$ with genus $g$ and 
and $n_1,\dots,n_k$ are the number of prongs
at the singularities and boundary components, then by
the Poincar\'e-Hopf theorem
\begin{eqnarray}\label{chi-eqn}
\sum_{i=1}^k (n_i - 2) = 4g-4.
\end{eqnarray}

For $(a,b) = (1,n)$, $n$ not divisible by $3$, we have two singularities
with number of prongs given by:
\begin{eqnarray*}
\psi_n(\gamma_1) &=& r_1 + 3n\\
\psi_n(\gamma_2) &=& 3s_2 + n.
\end{eqnarray*}
Plugging into (\ref{chi-eqn}) gives
$$
r_1 + 3s_2=0.
$$
The mapping class $\phi_{(1,2)}$ is the unique
genus 2 pseudo-Anosov mapping class with dilatation equal to $\lambda_2$ 
\cite{CH08, LT09},  and has one 6-pronged singularity \cite{HK:braidbounds}.  
Thus, $r_1= s_2=0$ and
$$
\gamma_1 = \ell_1 = 3\mu_2
$$
and
$$
\gamma_2 = \mu_2.
$$
The claim follows.
\qed

\begin{corollary}\label{prongs-cor} The map $\phi_{(a,b)}$ has singularities with
number of prongs (or {\it prong-type}) given by
$$
\left \{
\begin{array}{ll}
(3b,b) &\qquad\mbox{if $\gcd(3,ab) = 1$}\\
(3b,b/3,b/3,b/3) &\qquad\mbox{if $\gcd(3,b) = 3$}\\
(b,b,b,b) &\qquad\mbox{if $\gcd(3,a)= 3$}
\end{array}
\right .
$$
\end{corollary}

\begin{corollary}\label{orientability2-cor} If $b$ is odd, then $\phi_{(a,b)}$ is
not orientable.
\end{corollary}

\begin{corollary}\label{poles-cor} For $(a,b) \in \primPsi$, $\phi_{(a,b)}$ is 1-pronged
at one or more boundary components of $S_{(a,b)}$ if and only if 
$(a,b) \in \{(0,1),(\pm 1,3),(\pm 2,3)\}$.
\end{corollary}

\begin{corollary}\label{closures-cor}  If $(a,b)\not\in \{(0,1),(\pm 1,3),(\pm 2,3)\}$, then
$\phi_{(a,b)}$ extends to the closure of $S_{(a,b)}$ over the boundary components to 
a mapping class $\overline \phi_{(a,b)}$ with the same dilatation as $\phi_{(a,b)}$.
\end{corollary}

\begin{proposition} Table 2 below  describes the pairs $(a,b) \in \primPsi$ that give
rise to an orientable (or non-orientable) genus $g$ pseudo-Anosov mapping class.
(Here $g \ge 4$.)

\bigskip
\begin{center}
\begin{tabular}{|c|l|l|}
\hline
$g\ (\mod \ 6) $& orientable  & non-orientable\\
\hline
0 & no example & $b = g+1$, $a = 0 (\mod\ 3)$\\
\hline
1 & $b=g+1$, $a = 3(\mod \ 6)$ & $b=g$, $a=1,2(\mod\ 3)$\\
\hline
2 & $b = g$, $a = 1,5(\mod \ 6)$  & $b = g+1$, $a = 1,2(\mod \ 3)$\\
\hline
3 & $b = g+1$, $a= 3(\mod \ 6)$ & no example\\
\hline
4 & $b = g$, $a = 1,5(\mod \ 6)$ & $b = g+1$, $a = 0 (\mod \ 3)$\\
\hline
5 & $b = g+1$, $a = 1,5(\mod \ 6)$ & $b = g$, $a = 1,2(\mod \ 3) $\\
\hline
\end{tabular}
\bigskip

{Table 2:  Fibrations of $M$  according to genus.}
\end{center}
\end{proposition}

\section{Minimal dilatations for the fibered face.}\label{growth-section}

Let $\primPsi$ be the primitive elements of the
fibered cone discussed in Section~\ref{example-section}.
Let 
\begin{eqnarray*}
d_g &= &\min\{ \lambda(\psi) :  \psi \in \primPsi, \mbox{genus of $\psi$ is $g$}\}, \mbox{and}\\
d_g^+ &=& \min\{ \lambda(\psi) :  \psi \in \primPsi, \mbox{genus of $\psi$ is $g$, the
monodromy of $\psi$ is orientable}\}.
\end{eqnarray*}
In this section, we finish the proofs of Theorem~\ref{bounds-thm} and 
Theorem~\ref{orientable-bounds-thm} and their consequences
by determining $d_g$ and $d_g^+$.

\begin{proposition}\label{Fried-cor} Let $(a,b) \in \primPsi$.  
Then
$$
\lambda_{(a,b)} < \lambda_{(a',b')}
$$
if either
\begin{enumerate}
\item $|a|<|a'|$ and $|b|=|b'|$; or
\item $|a|=|a'|$ and $|b| > |b'|$.
\end{enumerate}
\end{proposition}

\noindent
{\bf Proof.} 
One compares the slopes of
rays from the origin to $(a,b)$ and $(a',b')$.  The claim follows from Theorem~\ref{Fried-thm}.
\qed


\begin{proposition}\label{compare-b-prop}  
For  $b\ge 3$, we have
$$
\lambda_{(1,b)} \ge \lambda_{(3,b+1)},
$$
with equality when $b=3$.
\end{proposition}

\noindent
{\bf Proof.}
Let $\lambda=\lambda_{(3,b+1)}$.
We will show that $LT_{(1,b)}(\lambda) < 0$.
Multiplying by $\lambda^2$ and using the fact that $LT_{(3,b+1)}(\lambda) = 0$ gives
\begin{eqnarray*}
\lambda^2 LT_{(1,b)}(\lambda) &=& \lambda^2LT_{(1,b)} (\lambda) - LT_{(3,b+1)}(\lambda)\\
&=&\lambda^{b+4} - \lambda^{b+3} - \lambda^{b+2} + \lambda^{b-2} + \lambda^2 - 1\\
&=&(\lambda-1)(\lambda^{b+3} - \lambda^{b-2}(\lambda^3 + \lambda^2 + \lambda + 1) + \lambda + 1)\\
&=&(\lambda-1)\lambda^{b-2}[\lambda^5-\lambda^3 - \lambda^2 -\lambda-1 + \lambda^{2-b}(\lambda+1)].
\end{eqnarray*}
Thus, it is enough to show that
$$
\lambda^5 - \lambda^3 - \lambda^2 -\lambda - 1 + \lambda^{2-b}(\lambda+1) < 0.
$$
Let $C$ be the quantity on the left side of this inequality. 
Since $\lambda > 1$ and $b \ge 3$, we have
$$
C < \lambda^5 - \lambda^3 - \lambda^2 = \lambda^2(\lambda^3 - \lambda - 1).
$$
One can check that the right hand side is negative for 
$$1 < \lambda < 1.3. $$
By Proposition~\ref{Fried-cor}, $\lambda$ decreases as $b$ increases.  A check shows that
$$1<\lambda_{(3,5)} < 1.3,$$
 and hence  $C<0$ for $b \ge 4$.  For $b=3$, one checks directly
that 
$$
\lambda_{(1,3)} = \lambda_{(3,4)}.
$$
\qed

\begin{remark}  The mapping class $\phi_{(1,3)}$ is defined on a genus $2$ surface
with four boundary components, with prong-type (3,1,1,1) and is not orientable.  The mapping class
$\phi_{(3,4)}$ is defined on a genus $3$ surface with prong-type (4,4,4,4) and is orientable.
By Proposition~\ref{compare-b-prop} these two examples have the same dilatation.
\end{remark}

Putting together Proposition~\ref{Fried-cor} and Proposition~\ref{compare-b-prop}, we have the following.

\begin{proposition}\label{monotonicity-prop} The  sequences $\lambda_{(1,b)}$ and $\lambda_{(3,b)}$ 
satisfy:
$$
\lambda_{(1,b)} > \lambda_{(3,b+1)} > \lambda_{(1,b+1)}.
$$
\end{proposition}

\begin{proposition} For $n \ge 2$, 
$$
\lim_{n \rightarrow \infty} (\lambda_{(a,n)})^n = \frac{3 + \sqrt{5}}{2},
$$
for any fixed $a$.
 \end{proposition}
 
 \noindent
 {Proof.} The rays through the lattice points $(a,n) \in \Lambda_M$ on the fibered
 face of $\psi$ converge to the ray through $(0,1)$.
 \qed
 
 \begin{corollary}\label{limit-d-cor}  For the minimal dilatations $d_g$ and $d_g^+$ 
 that are realized on $M$, we have
$$
\lim_{g \rightarrow \infty}\  (d_g)^g =  \frac{3+\sqrt{5}}{2},
$$
and
$$
\lim_{\begin{array}{c}\text{\tiny $g \rightarrow \infty$}\\ \text{\tiny $g \neq 0(\mod\ 6)$}\end{array}} (d_g^+)^g =  \frac{3+\sqrt{5}}{2}.
$$
\end{corollary}

\begin{figure}[htbp]
\begin{center}
\includegraphics[height=10cm]{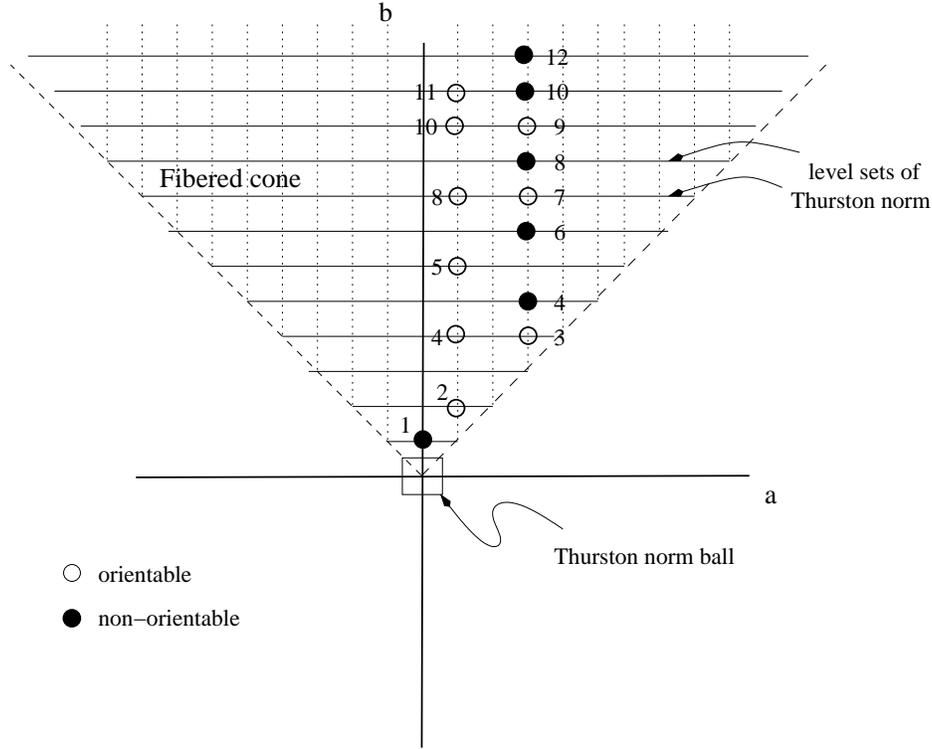}
\caption{Minima for $d_g$ and $d_g^+$ in genus $g =1,\dots,12$.}
\label{grid-fig}
\end{center}
\end{figure}
 
\begin{proposition}\label{bestbounds-prop} The following table describes the pairs $(a,b) \in  \primPsi$ that give
rise to the minima $d_g$ and $d_g^+$ realized on $M$. 

\bigskip
 \begin{center}
\begin{tabular}{|c|l|l|}
\hline
$g\ \mod \ 6 $& $\lambda(\phi_{(a,b)} )=d_g^+$, $\phi_{(a,b)}$ orientable  & $\lambda(\phi_{(a,b)})=d_g$\\
\hline
0 & no example &($3,g+1$)\\
\hline
1 & ($3,g+1$) & ($3,g+1$)\\
\hline
2 & ($1,g$) & ($1,g+1$)\\
\hline
3 & ($3,g+1$) &($3,g+1$)\\
\hline
4 & ($1,g$) & ($3,g+1$)\\
\hline
5 & ($1,g+1$) & ($1,g+1$)\\
\hline
\end{tabular}

\nopagebreak
\bigskip
Table 3: Pairs $(a,b)$ giving smallest dilatations for $\phi \in \Phi(\Msb)$.
\end{center}
\end{proposition}

Proposition~\ref{bestbounds-prop} and Corollary~\ref{closures-cor} complete the proofs
of Theorem~\ref{bounds-thm} and Theorem~\ref{orientable-bounds-thm}.  A pictorial view
of how the elements of $ \Psi$ giving the least dilatations for each genus up to 12 lie
on a fibered cone of $M$ is shown in 
Figure~\ref{grid-fig}.
\medskip

Putting together the results of this paper with those in \cite{AD10, KT10, LT09}, we see that for genus $g = 2,3,4,5,7,$ and $8$, 
$$
\delta_g^+ = \lambda_{(a,b)}
$$
where 
$$
(a,b) = 
\left \{
\begin{array}{lr}
(1,g) & \mbox{if $g = 2,3,4,$ or $8$}\\
(1,g+1) & \mbox{if $g = 5$}\\
(2,g+2) & \mbox{if $g = 7$}\\
\end{array}
\right .
$$
and
$$
\delta_6^+ \ge \lambda_{(1,6)}.
$$

These results suggest the following generalization to Question~\ref{LT-ques}.

\begin{question} For every $g \ge 2$,
is it true that 
$$
\delta_g^+ = \lambda_{(a,b)}
$$
for some $a,b$ with $b \ge g \ge a \ge 1$?
\end{question}

 \bibliographystyle{math}
 \bibliography{math}
\bigskip

\verse{
Eriko Hironaka\\
Department of Mathematics\\
Florida State University\\
Tallahassee, FL 32306-4510\\
U.S.A.
}

 \end{document}